\documentclass{amsart}

\newtheorem{theorem}{Theorem}[section]
\newtheorem{lemma}[theorem]{Lemma}

\theoremstyle{definition}
\newtheorem{definition}[theorem]{Definition}

\newtheorem{proposition}[theorem]{Proposition}
\newtheorem{corollary}[theorem]{Corollary}

\theoremstyle{remark}
\newtheorem{remark}[theorem]{Remark}

\numberwithin{equation}{section}

\newcommand {\sg} {\mathrm{sg}}

\begin{document}

\title{Polynomial Detection of Matrix Subalgebras}
%\date{\today}   

\author{Daniel Birmajer}
\address{ Department of Mathematics $\&$ Computer Science,
Nazareth College, 4245 East Avenue, Rochester, NY 14618}
%\curraddr{Department of Mathematics and Statistics,
%Case Western Reserve University, Cleveland, Ohio 43403}
\email{abirmaj6@naz.edu}

%    General info

%\date{January 1, 1994 and, in revised form, June 22, 1994.}

%\dedicatory{This paper is dedicated to our authors.}


\begin{abstract}
The double Capelli polynomial of total degree $2t$ is
\begin{equation*}
\sum \left\{ (\mathrm{sg}\, \sigma\tau)
  x_{\sigma(1)}y_{\tau(1)}x_{\sigma(2)}y_{\tau(2)}\cdots
  x_{\sigma(t)}y_{\tau(t)} 
|\; \sigma,\, \tau \in S_t\right\}.
\end{equation*}
It was proved by Giambruno-Sehgal and Chang that the double Capelli polynomial of total degree $4n$ is a polynomial
identity for $M_n(F)$. (Here, $F$ is a field and  $M_n(F)$ is the algebra of
$n \times n$ matrices over $F$).
Using a strengthened version of this result obtained by Domokos, we
show that the double Capelli polynomial of 
total degree $4n-2$ is a polynomial identity for any proper $F$-subalgebra of $M_n(F)$.
Subsequently, we present a similar result for nonsplit extensions of
full matrix algebras.
\end{abstract}
\maketitle
%%%%%%%%%%%%%%%%%%%%%%%%%%%%%%%%%%%%%%%%%%%%%%%%%%%%%%%%%%%%%%%%%%%%%%%%%%%%
\section{Introduction}\label{S:notation}
The double Capelli polynomial of total degree $2t$ is
\begin{equation*}
\sum \left\{ (\mathrm{sg}\, \sigma\tau)
  x_{\sigma(1)}y_{\tau(1)}x_{\sigma(2)}y_{\tau(2)}\cdots
  x_{\sigma(t)}y_{\tau(t)} 
|\; \sigma,\, \tau \in S_t\right\}.
\end{equation*}
In this paper we show that the double Capelli polynomial of degree $4n-2$
is a polynomial identity for any proper subalgebra of $M_n(F)$.
Subsequently, we present a polynomial test for nonsplit extensions
of full matrix algebras.

To begin, let $F$ be a field,  $M_n(F)$ the algebra of $n \times n$ matrices
over $F$, and 
$F\left\{X\right\}=F\left\{x_1, x_2, \dots\right\}$ the free associative algebra over $F$ in countably many variables. 
Sometimes we will use other variables $x$, $y$, $z$, $x_i$, $y_i$ for notation simplicity. 
A nonzero polynomial $f(x_1,\dots x_m) \in F\left\{X\right\}$ is a
\emph{polynomial identity} for an $F$-algebra $R$ if 
$f(r_1, \dots ,r_m)=0$ for all $r_1, \dots , r_m \in R$. A
\emph{$T$-ideal} is an ideal of $F\left\{X\right\}$ 
which is closed under endomorphisms of $F\left\{X\right\}$.
If $f_1,\dots, f_t$ are polynomial identities for $R$,
so is every polynomial $f$ in the $T$-ideal generated by $f_1, \dots ,f_t$. 
In this case we say that the identity $f=0$ in
$R$ is a \emph{consequence} of the identities $f_i=0,$ for $1 \le i \le t$.  

The standard polynomial of degree $t$ is
\begin{equation*}
s_t(x_1, \dots, x_t) = \sum_{\sigma \in S_t} (\sg\sigma) 
x_{\sigma(1)}x_{\sigma(2)} \cdots x_{\sigma(t)},
\end{equation*} 
where $S_t$ is the symmetric group on $\{1, \dots, t\}$ and $(\sg\sigma)$ is the sign of 
the permutation $\sigma \in S_t$.
The standard polynomial $s_t$ is homogeneous of degree $t$, multilinear and alternating.

The Amitsur-Levitzki theorem asserts that $M_n(F)$ satisfies any standard
polynomial of degree $2n$ or higher. Moreover, if $M_n(F)$ satisfies a polynomial of degree
$2n$, then the polynomial is a scalar multiple of $s_{2n}$ (cf. \cite{AL50}).
The \emph{Capelli polynomials} are 
\begin{equation*}
\begin{split}
c_{2t-1}\left(x_1, \dots ,x_t , y_1, \dots , y_{t-1}\right) &= \sum_{\sigma \in S_t} (\sg \sigma) 
x_{\sigma(1)}y_{1}x_{\sigma(2)}y_{2} \cdots x_{\sigma(t-1)}y_{t-1} x_{\sigma(t)},\\
\intertext{and}
c_{2t}\left(x_1, \dots ,x_t , y_1, \dots , y_{t}\right)
&=c_{2t-1}\left(x_1, \dots ,x_t , y_1, \dots , y_{t-1}\right)y_t.
\end{split}
\end{equation*}
These polynomials were introduced by Razmyslov in \cite{Ra74}.
The polynomials $c_{2t-1}$ and $c_{2t}$ are multilinear and alternating as a function of  $x_1, \dots , x_t$.
It is clear by a dimension argument that $c_{2n^2}$ is a PI for any proper
F-subalgebra of $M_n(F)$. On the other hand, $c_{2n^2}$ is not a PI for
$M_n(F)$. To see this, evaluate \\ $c_{2n^2}\left( x_1, \dots , x_{n^2}, y_1, \dots y_{n^2}\right)$ with
\begin{equation*}
\begin{split}
\left( x_1, x_2, \dots , x_{n}, x_{n+1}, \dots x_{n^2-1}, x_{n^2}\right)
&=\left(e_{11}, e_{12}, \dots , e_{1n}, e_{21}, \dots e_{n(n-1)}, e_{nn}\right),\\
\left( y_1, \dots , y_{n}, \dots y_{n^2-1}, y_{n^2} \right)
&=\left(e_{11}, \dots , e_{n2}, \dots e_{(n-1)n}, e_{n1} \right).
\end{split}
\end{equation*}  
where the $e_{ij}$ are the standard matrix units, $y_1=e_{11}$,
$y_{n^2}=e_{n1}$, and $y_2, \dots y_{n^2-1}$ are the 
unique choices of matrix units such that the monomial with 
$\sigma =1$ is nonzero, so $c_{2n^2}$ takes on the value $e_{11} \ne 0.$ 
Based on this example, we introduce the following definition:
\begin{definition}
We will say that a multilinear polynomial $f(x_1, \dots , x_t) \in F\{X\}$ is
a \emph{polynomial test} for an $F$-algebra $R$ if it is 
not a polynomial identity for $R$ but it is an identity for every proper $F$-subalgebra of $R$.   
\end{definition}
Thus, the Capelli polynomial of total degree $2n^2$ is a polynomial test for $M_n(F)$.
Moreover, central polynomials for $M_n(F)$  are polynomial tests for $M_n(F)$ (see \cite{Fo72}).
In \cite{Bi03}, it is proved that the standard polynomial of degree $2n-2$
is a polynomial test for the subalgebra of upper triangular matrices of $M_n(F)$.
The \emph{double Capelli polynomials} are
\begin{equation*}
\begin{split}
h_{2t-1}(x_1, \dots ,x_t , y_1,& \dots ,y_{t-1}) \\
&= \sum_{ \sigma \in S_t, \tau \in S_{t-1}} (\sg \sigma\tau) x_{\sigma(1)}y_{\tau(1)}x_{\sigma(2)}y_{\tau(2)} \cdots
x_{\sigma(t-1)}y_{\tau(t-1)}  x_{\sigma(t)}, 
\end{split}
\end{equation*}
and
\begin{equation*}
\begin{split}
h_{2t}(x_1, \dots ,x_t , y_1,& \dots , y_{t}) \\
&=\sum_{\sigma, \tau \in S_t} (\sg \sigma\tau)
x_{\sigma(1)}y_{\tau(1)}x_{\sigma(2)}y_{\tau(2)} 
\cdots x_{\sigma(t-1)}y_{\tau(t-1)}
x_{\sigma(t)}y_{\tau(t)}.
\end{split}
\end{equation*} 
Note that $h_{2t-1}$ and $h_{2t}$ are multilinear and alternate in the $x_i$ and also in the $y_j$.

Formanek pointed out that $h_{4n-2}$ is not a polynomial identity for $M_n(F)$
and asked for the least integer $m$ such that $h_m$ is a polynomial identity for $M_n(F)$.
Chang \cite{Ch88} has proved that the double Capelli polynomial
$h_{2t}$ is a consequence of the 
standard polynomial $s_t$.
A different proof that $h_{4n}=0$ is a polynomial identity for $M_n(F)$, that
uses a variation of Rosset's proof of the Amitsur-Levitzki theorem \cite{Ro76}, was given by
Giambruno-Sehgal in \cite{GS89}.
An elegant one-line proof of Domokos is given in \cite{Do93},\,\emph{\rm Example 2.2,
p. 917}.

In  \cite{Do95}, Domokos obtained a generalization of  Chang's theorem.
Since it is important in these notes, the precise statement of
Domokos's theorem is included below.  

Let $x_1,\dots,x_d,y_1,\dots,y_m$ be noncommuting variables over
$F$, and let $w_1,\dots,w_u$ be monomials in $y_1,\dots,y_m$ such
that $w_1, \dots ,w_u$ is a reordering of $y_1,\dots ,y_m$. 
For a subset $\Pi \subseteq S_d$ and a monomial partition
$\{w_1, \dots , w_u\}$ of the set of variables $Y$ we put
\begin{multline*}
f_\Pi(x_1,\dots ,x_d, y_1, \dots ,y_m |w_1, \dots ,w_u) =\\
 \sum (\sg \,\mu)x_{\pi(1)}\cdots x_{\pi(d_1)}w_{\rho(1)}
x_{\pi(d_1+1)}\cdots x_{\pi(d_1+d_2)}w_{\rho(2)}\cdots\\
\cdots w_{\rho(s)}x_{\pi(d_1+ \cdots +d_s+1)}\cdots x_{\pi(d_1+
\cdots +d_{s+1})},
\end{multline*}
where the summation runs  over all $\pi \in \Pi, \rho \in S_u$,  $d_i
\ge 1 $ for $i=1,\dots ,s+1$ such
that $d_1+\cdots d_{s+1}=d$ and  $\sg \,\mu$ is $\pm 1$ according to the
parity of the permutation of the ``underlying'' variables $x_1, \dots , x_d,
y_1, \dots ,y_m$ in the corresponding term.
\begin{theorem}\label{T:Domokos}
\cite{Do95}
The polynomial $f_{S_d}(x_1,\dots ,x_d, y_1, \dots , y_m | w_1, \dots
, w_u)$ is contained in the $T$-ideal generated by the standard
polynomial $s_d$.
\end{theorem} 
\begin{corollary}
\cite{Do95}
We have the strengthened version of the result of \cite{Ch88}
and \cite{GS89} we mentioned above:
\begin{equation*}
\sum_{ \sigma \in S_{2n}, \tau \in S_{2n-1}} (\sg \sigma\tau) x_{\sigma(1)}y_{\tau(1)} \cdots
y_{\tau(2n-1)} x_{\sigma(2n)}=0, 
\end{equation*} 
is a polynomial identity for $M_n(F)$, moreover, it is a consequence
of the standard identity $s_{2n}=0$.
\end{corollary}
To see that $h_{4n-2}$ is not a polynomial identity for $M_n(F)$, consider the substitution (double staircase)
\begin{multline*}
x_1=e_{11}, y_1=e_{12}, x_2=e_{22}, y_2=e_{23}, \dots , x_n=e_{nn}\\
y_n=e_{nn}, x_{n+1}=e_{n(n-1)}, y_{n+1}=e_{(n-1)(n-1)}, \dots , x_{2n-1}=e_{21}, y_{2n-1}=e_{11}
\end{multline*}
where the $e_{ij}$ are the standard matrix units. 
The only nonzero monomials in  $h_{4n-2}(x_i, y_i)$ are the $2n-1$ even cyclic
permutations of $x_1y_1 \dots x_{2n-1}y_{2n-1}$, 
and they all have positive sign. Thus
\begin{equation*}
h_{4n-2}(x_1, \dots ,x_{2n-1} , y_1, \dots , y_{2n-1})=2I-e_{11}.
\end{equation*}

We finish this section with two useful properties of the double 
Capelli polynomials.   
\begin{proposition}\label{P:Section1properties}
\begin{enumerate}
\item[(a)] $h_{q+r}$ is a linear combination, with coefficients being
$1$ or $-1$ of evaluations of $h_{q}\, h_{r}$.
\item[(b)] The polynomial $h_t$ is a consequence of the identity $h_s$
for any $t \ge s$.
\end{enumerate}
\end{proposition}
\begin{proof}
To prove (a) we show an explicit formula, where for simplicity we
consider the following statement:
$h_{2(q+r)-2}$ is a linear combination  with coefficients being
$1$ or $-1$ of evaluations of $h_{2q-1}\, h_{2r-1}$.
Let $t=q+r-1$
We  partition the set of permutations $S_{t}$ by defining the equivalence relation 
$\sigma_1 \sim_{q} \sigma_2$ if the images of the interval $[1,q]$
under $\sigma_1$ and $\sigma_2$ are the same set. Similarly, We
partition the set of permutations $S_{t}$ 
by defining the equivalence relation 
$\tau_1 \sim_{r} \tau_2$ if the images of the interval $[1,q-1]$
under $\tau_1$ and $\tau_2$ are the same set. Then we have
\begin{multline*}
h_{2t}(x_1, \dots , x_{t}, y_1, \dots , y_{t})=\\
\sum_{\substack{\bar\sigma \in S_{t}/\sim_q\\\bar\tau \in S_{t}/\sim_r}}
(\sg\sigma\tau)\, h_{2q-1}(x_{\sigma(1)}, \dots , x_{\sigma(q)},
y_{\tau(1)}, \dots , y_{\tau(q-1)})\\
h_{2r-1}(y_{\tau(1)}, \dots , y_{\tau(t)},x_{\sigma(q+1}), \dots , x_{\sigma(t)}).
\end{multline*}
The assertion in (b) follows immediately from (a).
\end{proof}
%%%%%%%%%%%%%%%%%%%%%%%%%%%%%%%%%%%%%%%%%%%%%%%%%%%%%%%%%%%%%%%%%%%%%%%%%%%%%%
\section{A polynomial test for the full matrix algebra}\label{S:Ptestfull}
The main goal of this section is to prove that $h_{4n-2}$ is a polynomial test for $M_n(F)$.
Before proceeding to the proof of this theorem we need some
preliminaries and notation (cf. \cite{Le02}). 
Let $\ell , m$ be positive integers such that $\ell + m = n$ and set 
\begin{equation*}
E_{(\ell,m)} (F) = 
\begin{bmatrix}
 M_\ell(F) & M_{\ell \times m}(F)\\
 0         &  M_m(F)
\end{bmatrix},
\end{equation*}
an $F$-subalgebra of $M_n(F)$.
\begin{enumerate}
\item[(i)] Associated to  $E_{(\ell,m)}(F)$ are canonical F-algebra homomorphisms 
$$\pi_\ell\colon  E_{(\ell ,m)} (F) \to M_\ell (F) \text{\quad  and
\quad} \pi_m\colon  E_{(\ell ,m)}(F) \rightarrow M_m(F).$$
Further identify $M_\ell(F)$ and $M_m(F)$ with 
\begin{equation*}
\begin{bmatrix}
 M_{\ell}(F)   & 0\\
 0             &  0
\end{bmatrix},
\begin{bmatrix}
 0         & 0\\
 0         &  M_{m}(F)
\end{bmatrix},
\end{equation*}
respectively.
\item[(ii)] Associated to a subalgebra $A$ of $E_{(\ell ,m)}(F)$ are
homomorphic image subalgebras $A_\ell$ and $A_m$ in $M_\ell (F)$ and 
$M_m (F)$ respectively.
\item[(iii)] Set 
\begin{equation*}
T_{(\ell,m)} (F) = 
\begin{bmatrix}
 0         & M_{\ell \times m}\\
 0         &  0
\end{bmatrix},
\end{equation*}
the Jacobson radical of $E_{(\ell ,m)}(F)$.
\item[(iv)] Recall that every $F$-algebra automorphism $\tau$ of $M_n(F)$
is \emph{inner} (i.e., there exists an invertible $Q$ in $M_n(F)$ 
such that $\tau(a)=QaQ^{-1}$ for all $a \in M_n(F)).$
We will say that two $F$-subalgebras $A, A^{\prime}$ of $M_n(F)$ are \emph{equivalent} provided there 
exists an automorphism $\tau$ of $M_n(F)$ such that $\tau(A) = A^{\prime}$.
\end{enumerate}
\begin{lemma}\label{L:h2t-2}
Let $A$ be a subalgebra of $E_{(\ell ,m)} (F)$ such that $A_\ell$ 
satisfies $h_{q}$  and $A_m$ satisfies $h_{r}$.
Then $A$ satisfies $h_{(q+r)}$.
\end{lemma}
\begin{proof}
The hypothesis that $A_\ell$ satisfies $h_q$ implies that the evaluation
of $h_q$ on $A$ consists of matrices of the form
\begin{equation*}
\begin{pmatrix}
 0         & *\\
 0         &  *
\end{pmatrix}.
\end{equation*}
Similarly, the hypothesis that $A_m$ satisfies $h_r$ implies that the
evaluation of $s_r$ on $A$ consists of matrices of the form
\begin{equation*}
\begin{pmatrix}
 *         & *\\
 0         &  0
\end{pmatrix}.
\end{equation*}
Thus $A$ satisfies $h_q\,h_r$. Since $h_{q+r}$ is a linear combination
of evaluations of $h_q\, h_r$, $A$ satisfies $h_{q+r}$.
\end{proof}
\begin{theorem}\label{T:polytest}
$h_{4n-2}$ is an identity for any proper subalgebra of $M_n(F)$.
\end{theorem}
\begin{proof}
 Let $A$ be a proper subalgebra of $M_n(F)$.
If $A$ is simple, then it is a a finite dimensional central simple algebra over its center $k$.
Let $K$ denote the algebraic closure of $k$;  then   
$A \otimes_k K$ is a simple $K$-algebra in a natural way (cf. \cite{Ro80},\,\S 1.8),  with 
$\dim_K \left(A \otimes_k K\right) = \dim_k(A)$.
Also, $A \otimes_k  K \cong M_t(K)$ for some $t \le n$.
Since  $A$ is a proper subalgebra of $M_n(F)$ it follows that $t < n$.
Hence, by the Amitsur-Levitzki theorem,  $A \otimes_k K$ satisfies $s_{2n-2}$.
Since $h_{4n-5}$ lies in the $T$-ideal generated by $s_{2n-2}$, we have
 that $h_{4n-5}(A)=0.$  
If $A$ is not simple, it can be embedded as $F$-algebra in
$E_{(\ell ,m)}(F)$ for some suitable positive integers $\ell$ and $m$ (with
$\ell + m = n$).
Since $h_{4\ell-1}$ and $h_{4m-1}$ are identities for $M_\ell(F)$ and
$M_m(F)$ respectively, we apply Lemma~\ref{L:h2t-2} to obtain that
$h_{4n-2}$ is an identity  for $A$.
\end{proof}
%%%%%%%%%%%%%%%%%%%%%%%%%%%%%%%%%%%%%%%%%%%%%%%%%%%%%%%%%%%%%%%%%
\section{A Polynomial test for $E_{(\ell,m)}$}\label{S:PtestElm}
In this section we show that the double Capelli polynomial $h_{4n-3}$ is
 a polynomial test for the subalgebra $E_{(\ell,m)}$
of $M_n(F)$ for any positive integers $\ell, m$ such that $\ell+m=n$.
\begin{proposition}\label{P:Elm}
$h_{4n-3}$ is an identity for every proper subalgebra of $E_{(\ell, m)}$.
\end{proposition}
\begin{proof}
We consider all possible proper subalgebras of $E_{(\ell,m)}(F)$.
Let first consider a subalgebra $A$ of $E_{(\ell, m)}$ such that 
$A_\ell$ is a proper subalgebra of $M_\ell(F)$. Then
$h_{4\ell-2}$ is an identity for $A_\ell$ as established in 
Theorem~\ref{T:polytest}, and $h_{4m-1}$ is an identity for  $M_m(F)$.
Thus, by Lemma~\ref{L:h2t-2}, $h_{4n-3}$ is an identity for  
\begin{equation*}
\begin{bmatrix}
 A_\ell & M_{\ell \times m}(F)\\
 0         &  M_m(F)
\end{bmatrix},
\end{equation*}
and consequently an identity for $A$.
Similarly,  $h_{4n-3}$ is an identity for every subalgebra of $E_{(\ell, m)}$
such that $A_m$ is a proper subalgebra of $M_m(F)$.
Clearly, $h_{4n-4}$ is an identity for the semisimple case 
\begin{equation*}
\begin{bmatrix}
 M_\ell(F) & 0\\
 0         &  M_m(F)
\end{bmatrix}.
\end{equation*}
It only remains to consider the case when the projections
$A\rightarrow A_\ell$ and $A\rightarrow A_m$ are equivalent
representations of $A$, which means that $A$ there is a fixed matrix
$T$ such that $TA{_\ell}T^{-1}=A_m$. It easily follows that in this case
$A$ is equivalent to the $F$-subalgebra of the form
\begin{equation*}
\left\{\begin{bmatrix}
 a         & c\\
 0         & a
\end{bmatrix}:
a, c \in M_\ell (F)\right\}.
\end{equation*}
In  \cite{Bi03}, Proposition~$2.5$, it is proved that the standard 
polynomial $s_{2\ell}$ is an identity for this algebra,
hence, $h_{2n-1}$  is an identity for $A$.
\end{proof}
\begin{remark}
The polynomial $h_{4n-3}$ is not an identity for $E_{(\ell, m)}$. 
For instance, if $n=3$ and $A= E_{(1, 2)}$, we have
\begin{equation*}
h_{9}\left(e_{11}, e_{11}, e_{12}, e_{22}, e_{22}, e_{23}, e_{33}, e_{33}, e_{32}\right)= 2e_{12}.
\end{equation*}
\end{remark}
\begin{remark}
The above ideas can be generalized to prove that the double Capelli polynomial
$h_{4n-t-1}$ is a polynomial test for the block upper triangular matrix
algebra
\begin{equation*}
\begin{pmatrix}
 M_{\ell_1}(F) &  &                 &    &   \\
   & M_{\ell_2}(F)&              & \text{\huge *} &  \\
   &    & \ddots   &   &   \\
    & \text{\Large 0}   &               &   & M_{\ell_t}(F)  \\
\end{pmatrix}.
\end{equation*}
\end{remark}
%%%%%%%%%%%%%%%%%%%%%%%%%%%%%%%%%%%%%%%%%%%%%%%%%%%%%%%%%%%%%%%%%%%%%%%%%%%%
\specialsection*{ACKNOWLEDGMENTS}
I am indebted to Professor Edward Letzter for his help and
guidance. I am also grateful to the referee for many helpful
comments and suggestions that have substantially improved  
these notes.
%%%%%%%%%%%%%%%%%%%%%%%%%%%%%%%%%%%%%%%%%%%%%%%%%%%%%%
\bibliographystyle{amsplain}

\end{document}